\documentclass[12pt]{article}
\usepackage{ulem}
\usepackage{appendix}
\usepackage{graphicx,amsmath,amsfonts,amssymb,color,mathrsfs, amsmath,amsthm}
\usepackage{epsfig}
\newcommand{\chuhao}{\fontsize{19pt}{\baselineskip}\selectfont}

\newcommand{\BOX}{\hfill $\Box$}

  \newtheorem{theorem}{Theorem}[section]
 \newtheorem{lemma}{Lemma}[section]
 \newtheorem{assumption}{Assumption}[section]

 \newtheorem{remark}{Remark}[section]

\usepackage{pifont}

\title{\bf\color{black} \chuhao{Exact tail asymptotics for fluid models driven by an M/M/c queue}}
\author{
Wendi Li\thanks{ School of Mathematics and Statistics, New Campus, Central
South University, Changsha, Hunan, 410083, P.R. China, E-mail:
liwendi@csu.edu.cn.
\newline \indent \hspace*{-2.5mm}
$^{**}$ Corresponding author.  School of Mathematics and Statistics, New Campus, Central
South University, Changsha, Hunan, 410083, P.R. China, E-mail:
liuyy@csu.edu.cn.
\newline \indent \hspace*{-2.5mm}
$^{***}$ School of Mathematics and Statistics, Carleton
University, 1125 Colonel By Drive, Ottawa, ON Canada K1S 5B6,
E-mail: zhao@math.carleton.ca. } ,
 Yuanyuan Liu$^{**}$ and Yiqiang Q. Zhao$^{***}$}
\date{July 8, 2018}
\begin{document}
\maketitle

\begin{abstract}
 In this paper, we investigate exact tail asymptotics for the stationary distribution of a fluid model driven by the $M/M/c$ queue, which is a two-dimensional queueing system with a discrete phase and a continuous level.
 We extend the kernel method to study tail asymptotics  of its stationary distribution, and a total of three types of exact tail asymptotics is identified from our study and reported in the paper. 
\end{abstract}

\section{Introduction}

Fluid flows have been widely used for modelling information flows in performance analysis of packet telecommunication systems. In this area, fluid queues with Markov-modulated input rates have played an important role in the recent development. In such a fluid model, the rate of information change is modulated according to a Markov process, evolving in the
background. Several references on the Markov-modulated fluid queues can be found in the literature, such as \cite{S98, N04, GLR13}. In these studies, the state space $N$ of the modulating Markov process is assumed to be finite, which put a restriction on applications. On the contrary, in this paper, we consider an infinite capacity fluid model driven by the $M/M/c$ queue, which is a specific birth-death process. First, we let $Z(t)$  be the state of a continuous-time Markov chain on a countable state space, the background process, at time $t$,  and let $X(t)$ be the fluid level in the queue at time $t$. Let $r_{Z(t)}$ denote the rate of change of the fluid level (or the \textit{net input rate}) at time $t$. Then, the dynamics of the fluid level  $X(t)$ are given by
\begin{equation*}
    \frac{dX(t)}{dt}=\left\{ \begin{array}{ll}
 r_{Z(t)}, & \mbox{if $X(t)> 0$ or $r_{Z(t)}\geq 0$,} \\
     0, & \mbox{if $X(t)= 0$   and $r_{Z(t)}< 0$.}
   \end{array}
   \right.
\end{equation*}

Fluid queues driven by infinite state Markov chains have been considered in the past by several authors. For instance, van Doorn and Scheinhardt, in \cite{DS97}, for the stationary distribution of the fluid queue driven by a birth-death process, they used orthogonal polynomials to solve an infinite system of differential equations under certain boundary conditions and provided the same integral expression obtained by Virtamo and Norros in \cite{VN94} and by Adan and Resing in \cite{AR96} in the case driven by the $M/M/1$ queue. Parthasarythy and Vijayashree, in \cite{PV02}, provided an expression, via an integral representation of Bessel functions, for the stationary distributions of the buffer occupancy and the buffer content, respectively, for a fluid queue driven by an $M/M/1$ queue. By using the Laplace transform, they obtained a system of differential equations, which led to a continued fraction and the solution of the stationary distribution. In \cite{BS02}, Barbot and Sericola provided an analytic expression for the stationary distribution of the fluid queue driven by an $M/M/1$ queue through the generating function technique. Analysis for the transient distribution of the fluid queue driven by an $M/M/1$ queue was reported by Sericola, Parthasarathy and Vijayashreehad in \cite{SPV05}.

Although methods for studying the stationary performance measures of the fluid queue driven by the $M/M/1$ queue are different in the above mentioned references, the expressions obtained through integral expressions are usually cumbersome and hard to be used directly for asymptotic properties for the stationary distribution. In this paper, we extend the kernel method to characterize exact tail asymptotics for the stationary distribution of the fluid model driven by an $M/M/c$ queue. The main contributions include:

\begin{enumerate}
\item An extension of the kernel method. The key idea of the kernel method was proposed by Knuth in \cite{K69} and further developed by Banderier \textit{et al.} in \cite{BBD02}. The method has been recently extended to study the exact tail behaviour for two-dimensional stochastic networks (or random walks with reflective boundaries) for both discrete and continuous random walks in the quarter plane, for example see Li and Zhao~\cite{LZ11} and Dai, Dawson and Zhao~\cite{DDZ15}, and references therein.
Compared to other methods, the kernel method, which has been successfully used to study tail behaviour of models with both level and background either discrete or continuous, does not require a determination or characterization of the entire unknown function in order to characterize the exact tail asymptotic properties in stationary distributions.
It is worthwhile to point out that the application of the kernel method to the fluid model driven by an $M/M/c$ queue is not straightforward and requires significant efforts, since in this case the level is continuous and the background is discrete.

\item An extension of the finding for the tail asymptotic behaviour in the stationary distribution of a fluid queue driven by a Markov chain. We show in Section~\ref{sec:5} that for the fluid model driven by an $M/M/c$ queue, a total of three types of exact tail asymptotic properties exists, in comparison with the finding by Govorun, Latouche and Remiche in \cite{GLR13}, in which they showed that for a fluid model driven by a finite state Markov chain there is only one type of tail asymptotic property. This is also an extension of the tail asymptotic behaviour in the stationary distribution of a fluid queue driven by an $M/M/1$ queue, since the tail asymptotic property given in Case (iii) of Theorems~\ref{the-tail-asy-3} and \ref{the-tail-asy-4} does not exist for the case of $c=1$.

\end{enumerate}

The rest of the paper is organized as follows: In Section~\ref{sec:2}, we describe the fluid model, define the notation and present the system of partial differential equations satisfied by the joint probability distribution function of the buffer level and of the state of the driving process. In this section, we also establish the fundamental equation based on the differential equations. Section~\ref{sec:3} is devoted to the discussion on properties of the branch points in the kernel equation and the analytic continuation of the unknown functions in terms of the kernel method. In Section~\ref{sec:4}, an asymptotic analysis of the two unknown functions is carried out. In Section~\ref{sec:5}, an characterization on exact tail asymptotic in the stationary distribution for the model is presented. We show that there exist three types of
tail asymptotic properties for the boundary, joint, and marginal distributions, respectively. These results are an extension of the single type behaviour found in \cite{GLR13} for the stationary density of the fluid queue driven by a finite state Markov chain. In Section~\ref{sec:6}, two special cases ($c=1$ and $c=2$) are further considered. Finally, in Section~\ref{sec:7}, we  make some concluding remarks to complete the paper.

\section{Model description and fundamental equation} \label{sec:2}

We consider the fluid model driven by an $M/M/c$ queueing system $\{Z(t), t\geq 0\}$, where $Z(t)$ denotes the queue length of the  M/M/c queue at time $t$. It is known that $Z(t)$  is a special birth-death process with the state space
$\mathbb{E}=\{0, 1, 2, \ldots \}$. Let $\lambda_{i}$ be the arrival rate and $\mu_{i}$ be the service rate in state $i$ for $\{Z(t), t\geq 0\}$. Then,
\[
  \lambda_{i}=\lambda> 0 \ \ \mbox{for any} \ \ i\geq 0,
\]
and with $\mu>0$,
\[
    \mu_{i}= \left \{ \begin{array}{ll}
i\mu, & \mbox{for $0 \leq i\leq c-1$}, \\
c \mu, & \mbox{for $i\geq c$}. \end{array} \right.
\]
Suppose $\lambda< c\mu$. Then the unique stationary distribution $\xi=(\xi_{i})_{i\in \mathbb{E}}$ of $\{Z(t)\}$ exists, which is given by
\begin{equation*}
 \xi_{i}= \left \{ \begin{array}{ll}
\xi_{0} \displaystyle \frac{\rho^{i}}{i!}, & \mbox{for $1 \leq i\leq c$}, \\
\xi_{c} \displaystyle \left (\frac{\rho}{c} \right )^{i-c}, & \mbox{for $i > c$}, \end{array} \right.
\end{equation*}
where  $\xi_{0}=(\sum_{i=0}^{c-1}\frac{\rho^{i}}{ i!}+\frac{\rho^{c}}{(c-1)!(c-\rho)})^{-1}$ and $\rho=\frac{\lambda}{\mu}$.

According to \cite{SPV05}, we may regard the fluid model driven by an $M/M/c$ queue $\{Z(t), t\geq 0\}$ as a fluid commodity, which is referred to as \textit{credit}. The credit accumulates in an infinite capacity buffer during the full busy period of $M/M/c$ queue (i.e. whenever a customer arrives and finds all servers busy) at a positive rate $r_{Z(t)}$, defined as  $r_{i}=r> 0$ for any $i\geq c$. The credit depletes the fluid during the partial busy period of $M/M/c$ queue (i.e. whenever an arriving customer finds less than $c$ customers in the queue) at a negative rate $r_{Z(t)}$. It is reasonable to assume that the negative rate $r_{i}$ increases in $i$. Without loss of generality, we assume that  the net input rate is $r_{i}= i-c$ for any $0\leq i\leq c-1$.

In order that the stationary distribution of $X(t)$ exists, we shall assume throughout the paper that
\[
 \sum_{i\in \mathbb{E}}\xi_{i}r_{i}<0,
\]
which is equivalent to
\begin{equation*}
   (r+1)\lambda < c\mu+ (c\mu-\lambda)\cdot\sum_{i=0}^{c-2}\frac{(c-i)\lambda^{i+1-c}\cdot(c-1)!}{\mu^{i+1-c}\cdot i!}.
\end{equation*}

Now, we denote
\[
 F_{i}(t,x)= P\{Z(t)=i, X(t)\leq x\}
\]
for any $t\geq 0,$ $x\geq 0$ and $i \in \mathbb{E}$.
It is well known (see e.g. \cite{DS97}) that the joint distribution $F_{i}(t,x)$ satisfies the following partial differential
equations:
\begin{eqnarray*}
 \nonumber \frac{\partial F_{0}(t,x)}{\partial t} &=& c\frac{\partial F_{0}(t,x)}{\partial x} -\lambda F_{0}(t,x) +\mu F_{1}(t,x), \\
 \label{equ-par-dif-1}  \frac{\partial F_{i}(t,x)}{\partial t} &=& (c-i)\frac{\partial F_{i}(t,x)}{\partial x} +\lambda F_{i-1}(t,x)-(\lambda+ i\mu)F_{i}(t,x) +(i+1)\mu F_{i+1}(t,x), \ \
 1\leq  i \leq c-1,\\
\label{equ-par-dif-2} \frac{\partial F_{i}(t,x)}{\partial t} &=& -r\frac{\partial F_{i}(t,x)}{\partial x} +\lambda F_{i-1}(t,x)-(\lambda+ c\mu)F_{i}(t,x) +c\mu F_{i+1}(t,x), \ \ i\geq c.
\end{eqnarray*}

Let $Z$ and $X$  be the stationary states of $Z(t)$ and $X(t)$ respectively. Then, the stationary distribution is given by
\[
  \Pi_{i}(x)=\lim_{t\rightarrow \infty}F_{i}(t,x) = P\{Z=i, X \leq x\}.
\]
Define $\pi_{i}(x)=\frac{\partial \Pi_{i}(x)}{\partial x}$ for any $x> 0$ and $\pi_{i}(0)=\lim_{x\rightarrow 0^{+}}\pi_{i}(x)$. From the above partial differential equations, we have the following equations:
 \begin{align}\label{equ-joi-dis-1}
 -c\pi_{0}(x) & = \mu\Pi_{1}(x)-\lambda\Pi_{0}(x), \\
\label{equ-joi-dis-2}
  -(c-i)\pi_{i}(x) &=
  \lambda\Pi_{i-1}(x)-(\lambda+ i\mu)\Pi_{i}(x) + (i+1)\mu\Pi_{i+1}(x), \quad \mbox{for $1 \leq i \leq c-1$}, \\
\label{equ-joi-dis-3}
  r \pi_{i}(x) & =
  \lambda\Pi_{i-1}(x)-(\lambda+ c\mu)\Pi_{i}(x) + c\mu\Pi_{i+1}(x), \quad \mbox{for $i\geq c$}.
\end{align}
The initial condition of (\ref{equ-joi-dis-1}), (\ref{equ-joi-dis-2}) and (\ref{equ-joi-dis-3}) is given by
\[
\Pi_{i}(0)=0, \ \ i\geq c.
\]
In addition, for any $i\in \mathbb{E}$, we have
\[
\Pi_{i}(\infty)=\lim_{x\rightarrow \infty}\Pi_{i}(x)= \xi_{i}.
\]

Let $\phi_{i}(\alpha)$ be the Laplace transform for $\pi_{i}(x)$, i.e.,
\[
\phi_{i}(\alpha)= \int_{0}^{\infty} \pi_{i}(x)e^{\alpha x}dx.
\]
For any $i\in \mathbb{E}$, we have
\begin{equation*}
 \int_{0}^{\infty}\Pi_{i}(x)e^{\alpha x}dx= \int_{0}^{\infty}\left [\Pi_{i}(0)+\int_{0^{+}}^{x}\pi_{i}(s)ds \right ]e^{\alpha x}dx=-\frac{1}{\alpha}\Pi_{i}(0)-\frac{1}{\alpha}\phi_{i}(\alpha).
\end{equation*}
Thus taking the Laplace transforms of $\Pi_{i}(x)$ and $\pi_{i}(x)$ in (\ref{equ-joi-dis-2}) and (\ref{equ-joi-dis-3}), we can get
\begin{align}
 -\phi_{c-1}(\alpha) 
\nonumber &= -\frac{\lambda}{\alpha} \left [\Pi_{c-2}(0)+\phi_{c-2}(\alpha) \right ]+
 \frac{\lambda+(c-1)\mu}{\alpha} \left [\Pi_{c-1}(0)+\phi_{c-1}(\alpha) \right ]-
 \frac{c\mu}{\alpha}[\Pi_{c}(0)+\phi_{c}(\alpha)],
\end{align}
and for any $i\geq c$
\begin{align}\nonumber
 r \phi_{i}(\alpha) 
 &= -\frac{\lambda}{\alpha} \left [\Pi_{i-1}(0)+\phi_{i-1}(\alpha) \right ]+
 \frac{\lambda+c\mu}{\alpha} \left [\Pi_{i}(0)+\phi_{i}(\alpha) \right ]-
 \frac{c\mu}{\alpha}[\Pi_{i+1}(0)+\phi_{i+1}(\alpha)].
\end{align}
It then follows that 
\begin{eqnarray}\label{equ-gene}
\nonumber && \sum_{i=c-1}^\infty [-\lambda z^{2}+(-\alpha r+\lambda+c\mu)z-c\mu] \phi_{i}(\alpha)z^i \\
\nonumber    &=& \lambda\phi_{c-2}(\alpha)z^{c}+[(\mu-\alpha-\alpha r)z-c\mu]\phi_{c-1}(\alpha) z^{c}+ \sum_{i=c-1}^\infty [\lambda z^{2}-(\lambda+ c\mu)z+c\mu]\Pi_{i}(0)z^i \\
 \nonumber    &&+\lambda\Pi_{c-2}(0)z^c+(\mu z-c\mu)\Pi_{c-1}(0)z^{c-1}.
\end{eqnarray}



Denote
\[
\psi(\alpha, z)=\sum_{i=c-1}^\infty \phi_{i}(\alpha) z^i,
\]
and
\[
\psi(z)=\sum_{i=c-1}^\infty \Pi_{i}(0)z^i.
\]
Then, we can obtain the following fundamental equation, which connects the bivariate unknown function $\psi(\alpha, z)$ to
 the univariate unknown functions $\phi_{c-2}(\alpha)$, $\phi_{c-1}(\alpha)$ and $\psi(z)$:
\begin{equation*}
 H(\alpha, z)\psi(\alpha, z)= \lambda z^{c}[\phi_{c-2}(\alpha)+ \Pi_{c-2}(0)]+ H_{1}(\alpha, z) \phi_{c-1}(\alpha)+  H_{2}(\alpha, z)\psi(z)+ H_{0}(\alpha, z)\Pi_{c-1}(0),
\end{equation*}
where
\begin{eqnarray*}
  &&  H(\alpha, z)= -\lambda z^{2}+(-\alpha r+\lambda +c \mu) z- c\mu, \\
  && H_{1}(\alpha, z)=(\mu-\alpha r -\alpha)z^{c}-c\mu z^{c-1}, \\
  &&  H_{2}(\alpha, z)= H_{2}(z)=\lambda z^{2}-\lambda z-c \mu z+ c\mu,\\
  && H_{0}(\alpha, z)=H_{0}(z)=\mu z^{c}- c\mu z^{c-1}.
\end{eqnarray*}

By establishing a relation between $\phi_{c-2}(\alpha)$ and $\phi_{c-1}(\alpha)$, we obtain the following result.

\begin{theorem}\label{the-phi-rel}
The fundamental equation can be rewritten as
 \begin{equation}\label{equ-funde}
   H(\alpha, z)\psi(\alpha, z)= \hat{H}_{1}(\alpha, z)\phi_{c-1}(\alpha)+  H_{2}(z)\psi(z)+ \hat{H}_{0}(\alpha, z),
 \end{equation}
 where
 \begin{eqnarray*}
  \hat{H}_{1}(\alpha, z) &=& \lambda z^{c}A_{c-2}(\alpha)+ H_{1}(\alpha, z), \\
    \hat{H}_{0}(\alpha, z) &=& H_{0}(z)\Pi_{c-1}(0)+\lambda z^{c}\Pi_{c-2}(0)+\lambda z^{c}\sum_{n=0}^{c-2} \left[k_{n}\lambda^{c-2-n}\prod_{m=n}^{c-2}\frac{A_{m}(\alpha)}{(m+1)\mu}\right],
 \end{eqnarray*}
 with  $k_{0}=\mu\Pi_{1}(0)-\lambda\Pi_{0}(0)$, 
\[
 k_{i}=\lambda\Pi_{i-1}(0)-(\lambda+i\mu)\Pi_{i}(0)+ (i+1)\mu \Pi_{i+1}(0),\ 1\leq i\leq c-2,
\]
and
\[
 A_{i}(\alpha)=\frac{(i+1)\mu}{\alpha+\lambda+i\mu-\lambda A_{i-1}(\alpha)},  \ 0\leq i\leq c-2, \ \ A_{-1}(\alpha)=0.
\]

\end{theorem}

\proof
Taking the Laplace transform for $\Pi_{i}(x)$ and $\pi_{i}(x)$ in  (\ref{equ-joi-dis-1}) and (\ref{equ-joi-dis-2}),
 leads to the following linear equations:
\begin{equation*}
    \left\{
   \aligned
    & (\alpha+\lambda)\phi_{0}(\alpha)-\mu \phi_{1}(\alpha)= k_{0},
    \\
    & -\lambda\phi_{0}(\alpha)+(\alpha+\lambda+\mu)\phi_{1}(\alpha)-2\mu \phi_{2}(\alpha)= k_{1},
    \\
    &  \vdots
    \\
    & -\lambda\phi_{c-3}(\alpha)+[\alpha+\lambda+(c-2)\mu]\phi_{c-2}(\alpha)= (c-1)\mu\phi_{c-1}(\alpha)+k_{c-2}.
   \endaligned
   \right.
\end{equation*}
Since $A_{0}'(\alpha)=\frac{-\mu}{(\alpha+\lambda)^{2}}< 0$, we assume that $A_{k-1}'(\alpha)< 0$ for any $\alpha\geq  0$, as the inductive hypothesis, to show
\[
 A_{k}'(\alpha)=\frac{-(k+1)\mu[1-\lambda A_{k-1}'(\alpha)]}{[\alpha+\lambda+ k\mu-\lambda A_{k-1}(\alpha)]^{2}}< 0.
\]
Thus, $A_{i}(\alpha)$ is a decreasing function about $\alpha$ for any $0 \leq i\leq c-2$. For any $\alpha > 0$ and $0 \leq i\leq c-2$, we can obtain that
\[
 A_{i}(\alpha)< A_{i}(0)=\frac{(i+1)\mu}{\lambda},
\]
which implies that $A_{i+1}(\alpha)=\frac{(i+2)\mu}{\alpha+\lambda+(i+1)\mu-\lambda A_{i}(\alpha)}> 0$.
Hence  $0 < A_{i}(\alpha) < \frac{(i+1)\mu}{\lambda}$ for any $0 \leq i\leq c-2$ and $\alpha > 0$.

From the linear equations and the definition of $A_{i}(\alpha)$, we have for any $0 \leq i\leq c-2$,
\begin{equation}\label{equ-phi-rel}
 \phi_{i}(\alpha)= \sum_{n=0}^{i}[k_{n}\lambda^{i-n}\prod_{m=n}^{i}\frac{A_{m}(\alpha)}{(m+1)\mu}]+ A_{i}(\alpha)\phi_{i+1}(\alpha).
\end{equation}
Specially, for the case $c=1$, we have $ \hat{H}_{1}(\alpha, z) =  H_{1}(\alpha, z)$ and $\hat{H}_{0}(\alpha, z) = H_{0}(z)\Pi_{0}(0)$.
Hence, the theorem is proved.
\BOX

\section{Kernel equation and branch points} \label{sec:3}

The tail asymptotic behaviour of the stationary distribution for the fluid queue relies on properties of the kernel function $H(\alpha, z)$, and the functions $ \hat{H}_{1}(\alpha, z)$ and $H_{2}( z)$. Now, we consider the kernel equation
\begin{equation*}
  H(\alpha, z)= 0,
\end{equation*}
which can be written as a quadratic form in $z$ as follows
\begin{equation}\label{equ-qua-form}
  H(\alpha, z) = az^{2}+ b(\alpha) z+ d =0,
\end{equation}
where $a=-\lambda$, $b(\alpha)=-\alpha r+\lambda +c\mu$ and $d=- c\mu$.

Let
\[
 \Delta(\alpha)=b^{2}(\alpha)-4ad
\]
be the discriminant of the quadratic form in  (\ref{equ-qua-form}).
In the complex plane $\mathbb{C}$,
for each $\alpha$, the two solutions to (\ref{equ-qua-form}) are given by
\begin{equation}\label{equ-Z}
  Z_{\pm}(\alpha)=\frac{-b(\alpha)\pm \sqrt{\Delta(\alpha)}}{2a}.
\end{equation}
When $\Delta(\alpha)=0$, $\alpha$ is called a branch point of $Z(\alpha)$.

Symmetrically, for each $z$, the solution to (\ref{equ-qua-form}) is  given by
\begin{equation}\label{equ-alpha-1}
 \alpha(z)=\frac{-\lambda z^{2}+ (\lambda+c\mu)z-c\mu}{zr}.
\end{equation}
Note that all functions and variables are treated as complex ones throughout the paper.
We have the following property on the branch points.
\begin{lemma}\label{lem-bran-point}
$\Delta(\alpha)$ has two positive zero points $\alpha_{1} = \frac{\left ( \sqrt{c\mu} - \sqrt{\lambda} \right )^2}{r}$ and $\alpha_{2} = \frac{\left ( \sqrt{c\mu} + \sqrt{\lambda} \right )^2}{r}$. 
Moreover, $\Delta(\alpha)> 0$  in  $(-\infty, \alpha_{1})\cup (\alpha_{2}, \infty)$ and $\Delta(\alpha)< 0$  in $(\alpha_{1}, \alpha_{2})$.
\end{lemma}

For convenience, define the cut plane $\widetilde{\mathbb{C}}_{\alpha}$ by
\[
  \widetilde{\mathbb{C}}_{\alpha}= \mathbb{C}_\alpha \setminus \{[\alpha_{1},\alpha_{2}]\}.
\]
In the cut plane $\widetilde{\mathbb{C}}_{\alpha}$, denote the two branches of $Z(\alpha)$ by $Z_{0}(\alpha)$ and $Z_{1}(\alpha)$,
where $Z_{0}(\alpha)$ is the one with the smaller modulus and $Z_{1}(\alpha)$ is the one with the larger modulus. Hence we have
\[
  Z_{0}(\alpha)=Z_{-}(\alpha) \ \mbox{and} \ Z_{1}(\alpha)=Z_{+}(\alpha) \  \mbox{if} \  \Re(\alpha)> \frac{\lambda+c\mu}{r},
\]
\[
  Z_{0}(\alpha)=Z_{+}(\alpha) \ \mbox{and} \ Z_{1}(\alpha)=Z_{-}(\alpha) \  \mbox{if} \ \Re(\alpha) \leq \frac{\lambda+c\mu}{r}.
\]

\begin{lemma}\label{lem-ana-H}
The functions  $Z_{0}(\alpha)$ and $Z_{1}(\alpha)$ are analytic in $\widetilde{\mathbb{C}}_{\alpha}$.
Similarly, $\alpha(z)$
is meromorphic in $\mathbb{C}_z$ and $\alpha(z)$ has two zero points and one pole.
\end{lemma}

\proof
We first give a proof to $Z_{0}(\alpha)$ and the proof to $Z_{1}(\alpha)$ can be given in the same fashion.
Let $\alpha= a+bi$ with $a, b\in \mathbb{R}$ and $\arg(\alpha)\in (-\pi, \pi]$, and write $\Delta(\alpha)=\Re(\Delta(\alpha)) + \Im(\Delta(\alpha)) i$.   We then have
\[
    \Re(\Delta(\alpha))=R(a, b)=(a^{2}-b^{2})r^{2}-2(\lambda+c\mu)ra+(\lambda-c\mu)^{2},
\]
and
\[
 \Im(\Delta(\alpha))=I(a, b)=2abr^{2}-2(\lambda+c\mu)br.
\]
%
Let $\Im(\Delta(\alpha))=0$, we obtain that $a= \frac{\lambda+c\mu}{r}$ or $b=0$.
For $b=0$, from Lemma \ref{lem-bran-point}, we know that $R(a, b)\leq 0$ and $I(a, b)=0$ along the curve $\mathcal{C}_{1}=\{\alpha= a+bi: \alpha_{1}\leq a\leq\alpha_{2}, b=0\}$. According to the property of the square root function, if we take $\mathcal{C}_{1}$ as a cut of $\sqrt{\Delta(\alpha)}$, then the function $Z_{0}(\alpha)$ cannot be analytic on the curve $\mathcal{C}_{1}$. Thus, we will consider the analytic property for $Z_{0}(\alpha)$ on the cut plane $\widetilde{\mathbb{C}}_{\alpha}=\mathbb{C}_\alpha \setminus \mathcal{C}_{1}$ in the following.

For $a= \frac{\lambda+c\mu}{r}$ and any $b\in \mathbb{R}$, we obtain that
\[
    R\left (\frac{\lambda+c\mu}{r}, b \right )= R\left (\frac{\lambda+c\mu}{r}, 0 \right )-b^{2}r^{2}< R \left (\frac{\lambda+c\mu}{r}, 0 \right )< 0.
\]
Therefore, along the curve $\mathcal{C}_{2}=\{\alpha= a+bi: a=\frac{\lambda+c\mu}{r}\}$, we have $R(a, b)\leq 0$ and $I(a, b)=0$, which implies that $\sqrt{\Delta(\alpha)}$ or ($-\sqrt{\Delta(\alpha)}$) cannot be analytic on $\mathcal{C}_{2}$. However, from the definition of $Z_{0}(\alpha)$, we have that the branch $Z_{0}(\alpha)=Z_{+}(\alpha)$ is analytic in the domain $\{\alpha\in \widetilde{\mathbb{C}}_{\alpha}: \Re(\alpha) < \frac{\lambda+c\mu}{r}\}$ 
 and $Z_{0}(\alpha)=Z_{-}(\alpha)$ is analytic in the complementary domain of
the closure of this set in $\widetilde{\mathbb{C}}_{\alpha}$. 
From the choice of the square root, we know that the function $Z_{0}(\alpha)$ is continuous on the curves $\mathcal{C}_{2}$, 
which separates the two above domains. Thus, by Morera's Theorem, we have that the function $Z_{0}(\alpha)$ is analytic in the cut plane $\widetilde{\mathbb{C}}_{\alpha}$.


From (\ref{equ-alpha-1}), $\alpha(z)$ is analytic in $\mathbb{C}_z$ except at the pole $z=0$,
which implies that $\alpha(z)$ is meromorphic in $\mathbb{C}_z$. It also follows from (\ref{equ-alpha-1}) that $\alpha(z)$ has two zero points.
\BOX

Based on Lemma \ref{lem-ana-H}, we have the analytic continuation of $\hat{H}_{1}(\alpha, Z_{0}(\alpha))$ and $H_{2}(z)$.
\begin{lemma}\label{lem-ana-H1-H2}
 The function $\hat{H}_{1}(\alpha, Z_{0}(\alpha))$ is analytic on $\widetilde{\mathbb{C}}_{\alpha}$ and  $H_{2}(z)$
is analytic on $\mathbb{C}_z$.
\end{lemma}

\proof
From Theorem \ref{the-phi-rel}, we have
\[
\hat{H}_{1}(\alpha, Z_{0}(\alpha))=[\lambda A_{c-2}(\alpha)+\mu-\alpha r -\alpha]Z_{0}(\alpha)^{c}-c\mu Z_{0}(\alpha)^{c-1}.
\]
The analytic property is immediate from Lemma \ref{lem-ana-H}.

From the definition of $H_{2}( z)$, we can easily get the assertion.
\BOX

\section{Asymptotic analysis of $\phi_{c-1}(\alpha)$ and $\psi(z)$} \label{sec:4}

In order to characterize the exact tail asymptotics for the stationary distribution $\Pi_{i}(x)$, we need to study the asymptotic property of the two unknown functions $\phi_{c-1}(\alpha)$ and $\psi(z)$ at their dominant singularities, respectively. There are three steps in the asymptotic analysis of $\phi_{c-1}(\alpha)$ and $\psi(z)$:
(i) analytic continuation of the functions $\phi_{c-1}(\alpha)$ and $\psi(z)$; (ii) singularity analysis of the functions
$\phi_{c-1}(\alpha)$ and $\psi(z)$; and (iii) applications of a Tauberian-like theorem. In this section, we give details of the first and second steps, and the detail of the third step will be given in Appendix A.

We first introduce the following lemma, which is a transformation of Pringsheim's theorem for a generating function (see, for example,  Dai and Miyazawa \cite{DM11}).
\begin{lemma}\label{lem-radiu-conv}
Let $g(x)=\int_{0}^{\infty}e^{xt}f(t)dt$ be the moment generating function with real variable $x$.  The convergence parameter
of $g(x)$ is given by
\[
  C_{p}(g)=\sup \{x\geq 0: g(x)< \infty\}.
\]
Then, the complex variable function $g(\alpha)$ is analytic on $\{\alpha\in \mathbb{C}: \Re(\alpha)< C_{p}(g)\}$.
\end{lemma}

Now, we provide detailed information about the extended generator for the fluid queue, which will be used later to investigate the analytic continuation of $\phi_{c-1}(\alpha)$. 
Instead of focusing on the case that the modulated process is an M/M/c queue, we will consider a general setting, whose background process is a general continuous-time Markov chain with an irreducible,  conservative and countable (finite or infinitely countable) generator $Q=(q_{ij})$.  We first recall some related definitions.  Let $\Phi_{t}$ be a continuous-time Markov process with a locally compact, separable metric space $X$ and transition function $P^{t}(i,j)$. We denote by $\mathcal{D}(\mathcal{A})$ the set of all functions $f$, for which there exists a measurable function $g$ such that
the process $C_{t}^{f}$, defined by
\[
  C_{t}^{f}= f(\Phi_{t})-f(\Phi_{0})-\int_{0}^{t}g(\Phi_{s})ds,
\]
is a local martingale. We write $\mathcal{A}f=g$ and call $\mathcal{A}$ the extended generator of the  process $\Phi$.

Consider a general fluid model $(X(t), Z(t))$.   Define its weakly infinitesimal generator $\mathcal{B}$ of the fluid queue by
\[
 \mathcal{B}g(x,i)=\lim_{t\rightarrow 0}\frac{E_{(x,i)}[g(X(t), Z(t))]-g(x,i)}{t}.
\]
Then, we  present the following lemma.
\begin{lemma}\label{lem-ext-gen}
Let $(X(t), Z(t))$ be the general  fluid queue with the generator $Q=(q_{ij})$, and let $g(x,i)$ be a function such that $g$ is partially differentiable about $x$
and for any $x\geq 0$,
\begin{equation*}
  \sum_{j\in \mathbb{E}}g(x,j)q_{ij}<\infty.
\end{equation*}
Moreover, we assume that $\sup_{i\in  \mathbb{E}}|r_{i}|<\infty$.

(i) For $x> 0,$ $i\in \mathbb{E}$ or $x=0,$ $i\in \mathbb{E^{+}}$, we have
\[
 \mathcal{B}g(x,i)= r_{i}\frac{d g(x,i)}{d x}+ \sum_{j\in \mathbb{E}}g(x,j)q_{ij},
\]
and for $x=0,$ $i\in \mathbb{E^{-}}\cup \mathbb{E^{\circ}}$, we have
\[
 \mathcal{B}g(0,i)= \sum_{j\in \mathbb{E}}g(0,j)q_{ij},
\]
where $\mathbb{E^{+}}=\{i\in \mathbb{E}| r_{i}> 0\}$, $\mathbb{E^{-}}=\{i\in \mathbb{E}| r_{i}< 0\}$ and $\mathbb{E^{\circ}}=\{i\in \mathbb{E}| r_{i}= 0\}$.

(ii) If the partial derivative $\frac{d g(x,i)}{d x}$ is continuous in $x$, then $f\in \mathcal{D}(\mathcal{A})$ and $\mathcal{A}f=\mathcal{B}f$.
\end{lemma}

\proof
This proof is similar to the proof of Lemma 3.1 in \cite{LL18} and we will omit the detail here. It is worth noting that
 the phase process in \cite{LL18} is a finite continuous-time Markov chain, which is different from the phase process $\{Z(t)\}$ in this paper. In order to
extend the result in \cite{LL18}, we need to impose the assumption that $\sup_{i\in  \mathbb{E}}|r_{i}|<\infty$.
\BOX

According to Lemma \ref{lem-ext-gen}, we can state the following lemma, which is crucial for the analytic continuation of $\phi_{c-1}(\alpha)$ and $\psi(z)$.


\begin{lemma}\label{lem-ana}
$\phi_{c-1}(\alpha)$ is analytic on $\{\alpha: \Re(\alpha) < \alpha^{\ast}\}$, where $\alpha^{\ast}=C_{p}(\phi_{c-1})> 0$, 
and $\psi(z)$ is analytic on the disk $\Gamma_{z^{\ast}}=\{z: |z|< z^{\ast}\}$, where $z^{\ast}=\frac{c\mu}{\lambda}$. Moreover,
the following equation is satisfied in the domain $D_{\alpha, z}=\{(\alpha, z): H(\alpha, z)=0 \ \ and \ \ \psi(\alpha, z)< \infty\}$:
 \begin{equation}\label{equ-ana-rel}
   \hat{H}_{1}(\alpha, z)\phi_{c-1}(\alpha)+  H_{2}( z)\psi(z)+ \hat{H}_{0}(\alpha, z)=0.
 \end{equation}
\end{lemma}

\proof
First, we prove that $\alpha^{\ast}> 0$.  It follows from Lemma \ref{lem-ext-gen} that the extended generator is given by
\[
 \mathcal{A}V(x, i) = -e^{\alpha x}z^{i} \left [\lambda+c\mu-\alpha r-\lambda z-\frac{c\mu}{z} \right ],
\]
for $x\geq 0,$ $i\geq c$, and
\[
 \mathcal{A}V(x, i) = -e^{\alpha x}z^{i} \left [(c-i)\alpha+\lambda+i\mu -\lambda z-\frac{i\mu}{z} \right ],
\]
for $x>0,$ $0\leq i\leq c-1$, and
\[
   \mathcal{A}V(0, i) = -z^{i} \left [\lambda+i\mu-\lambda z-\frac{i\mu}{z} \right ],
\]
for $0\leq i\leq c-1$.

In order to find some constant $s> 0$ such that
\[
 \mathcal{A}V(x, i) \leq -sV(x, i),
\]
for $x> 0,$ $i\geq 0$,
we need to choose appropriate $\alpha$ and $z$ such that for any $0 \leq i\leq c-1$,
\begin{equation*}
\left \{
 \begin{array}{cc}
 \lambda+c\mu-\alpha r-\lambda z-\frac{c\mu}{z}>0,
\\
 (c-i)\alpha+\lambda+i\mu-\lambda z-\frac{i\mu}{z}>0.
   \end{array}\right.
\end{equation*}
For $\alpha=\frac{c\mu-\lambda}{r}> 0$, we have
\[
 \frac{\lambda+c\mu-r\alpha-\sqrt{\Delta_{1}}}{2\lambda}< 1,\ \ \frac{\lambda+i\mu+(c-i)\alpha-\sqrt{\Delta_{2}}}{2\lambda}< 1,
\]
and
\[
 \frac{\lambda+c\mu-r\alpha+\sqrt{\Delta_{1}}}{2\lambda}> 1,\ \ \frac{\lambda+i\mu+(c-i)\alpha+\sqrt{\Delta_{2}}}{2\lambda}> 1,
\]
where $\Delta_{1}=(\lambda+c\mu-r\alpha)^{2}-4c\lambda\mu$ and $\Delta_{2}=[\lambda+i\mu+(c-i)\alpha]^{2}-4i\lambda\mu$.

%
%
Thus, there exists some $z \in B= B_{1}\cap B_{2}\cap (1, \infty)\neq \varnothing $ such that
\begin{equation}\label{equ-dri-con}
  \mathcal{A}V(x, i)\leq -s V(x, i)+ b I_{\mathbb{L}_{0}},
\end{equation}
 where
\begin{eqnarray*}
  B_{1} &=& \left (\frac{\lambda+c\mu-r\alpha-\sqrt{\Delta_{1}}}{2\lambda},
\frac{\lambda+c\mu-r\alpha+\sqrt{\Delta_{1}}}{2\lambda} \right ),  \\
  B_{2} &=& \left (\frac{\lambda+i\mu+(c-i)\alpha-\sqrt{\Delta_{2}}}{2\lambda},
\frac{\lambda+i\mu+(c-i)\alpha+\sqrt{\Delta_{2}}}{2\lambda} \right ),\\
 s &=& \min \left \{\lambda+c\mu-\alpha r-\lambda z-\frac{c\mu}{z}, \lambda+i\mu+\alpha-\lambda z-\frac{i\mu}{z} \right \}> 0,\\
 \mathbb{L}_{0} &=& \{(x,i)| x=0, 0 \leq i \leq c-1\}.
\end{eqnarray*}
Since the drift condition (\ref{equ-dri-con}) holds, from Theorem 7 in \cite{MT}, we know that
\[
\phi_{c-1} \left (\frac{c\mu-\lambda}{r}\right )z^{c-1}< \psi\left (\frac{c\mu-\lambda}{r}, z \right )<\sum_{i=0}^{\infty}\int_{0}^{\infty}\pi_{i}(x)V(x,i)dx <\infty.
\]
Thus, from Lemma \ref{lem-radiu-conv} we can obtain that $\alpha^{\ast}\geq \frac{c\mu-\lambda}{r}> 0$.

For $\psi(z) = \sum_{i=c-1}^{\infty}\Pi_{i}(0)z^{i} $, we have
\[
 \psi(z) \leq \Pi_{c-1}(0)z^{c-1}+\sum_{i=c}^{\infty}\xi_{i} z^{i} = \Pi_{c-1}(0)z^{c-1}+ \xi_{c}z^{c}\sum_{i=0}^{\infty}(\frac{\lambda}{c\mu})^{i} z^{i}.
\]
Since $\sum_{i=0}^{\infty}(\frac{\lambda}{c\mu})^{i} z^{i}$ is convergent in $|z|< \frac{c\mu}{\lambda}$, we have that $\psi(z)$ is analytical in the disk $\Gamma_{\frac{c\mu}{\lambda}}$.

Now we prove the second assertion. From the equation (\ref{equ-funde}), we can obtain that if both $\phi_{c-1}(\alpha)$ and $\psi(z)$ are finite, then
$\psi(\alpha, z)$ is finite as long as $H(\alpha, z)\neq 0$. 
Assume that $H(\alpha_{0}, z_{0})=0$ for some $\alpha_{0}> 0$ and $1< z_{0}< \frac{c\mu}{\lambda}$, and $\phi_{c-1}(\alpha_{0})<\infty$, $\psi(z_{0})<\infty$. Then, for small enough $\varepsilon> 0$, we have  $\psi(\alpha_{0}, z_{0})< \psi(\alpha_{0}, z_{0}+\varepsilon)< \infty$,  thus (\ref{equ-ana-rel}) holds for such pair $(\alpha_{0}, z_{0})$.
\BOX

\begin{remark}\label{rem-dri}
Actually, accoridng to \cite{DMT95}, we know that the equation (\ref{equ-dri-con}) implies the fluid model driven by an $M/M/c$ queue is
$V$-uniformly ergodic.
 \end{remark}

Now, we present another relationship between  $\phi_{c-1}(\alpha)$ and $\psi(z)$, and extend their analytic domains.
\begin{lemma}\label{lem-asy-ana-1}

(i) $\phi_{c-1}(\alpha)$ can be analytically continued to the domain $D_{\alpha}=\{\alpha\in \widetilde{\mathbb{C}}_{\alpha}:  \hat{H}_{1}(\alpha, Z_{0}(\alpha))\neq 0\}\cap\{\alpha\in \widetilde{\mathbb{C}}_{\alpha}: |Z_{0}(\alpha) |< \frac{c\mu}{\lambda}\}$, and
\begin{equation}\label{equ-asy-ana-1}
   \phi_{c-1}(\alpha)=-\frac{H_{2}( Z_{0}(\alpha))\psi(Z_{0}(\alpha))+ \hat{H}_{0}(\alpha, Z_{0}(\alpha))}{ \hat{H}_{1}(\alpha, Z_{0}(\alpha))}.
\end{equation}

(ii) $\psi(z)$ can be analytically continued to the domain $D_{z}=\{z\in \mathbb{C}: H_{2}(z)\neq 0\}\cap\{z\in \mathbb{C}: Re(\alpha(z)) < \alpha^{\ast}\}$ and
\begin{equation}\label{equ-asy-ana-2}
   \psi(z) =-\frac{\hat{H}_{1}(\alpha(z), z) \phi_{c-1}(\alpha(z))+\hat{H}_{0}(\alpha(z), z)}{ H_{2}(z)}.
\end{equation}
\end{lemma}

\proof
%
(i) For any $(\alpha, z)$ such that $H(\alpha, z)=0$ and $\psi(\alpha,z)<\infty$, we can get equation (\ref{equ-ana-rel}).
Using $z=Z_{0}(\alpha)$ leads to (\ref{equ-asy-ana-1}).
Then, from Lemma \ref{lem-ana}, we know that the right-hand side of the above equation is analytic except for the points such that $\hat{H}_{1}(\alpha, Z_{0}(\alpha))=0$ or $|Z_{0}(\alpha)|\geq \frac{c\mu}{\lambda}$. Hence, we get the assertion.

Similarly, we can prove assertion (ii).

\BOX

%
Based on the above arguments, we have the following lemma.
\begin{lemma}\label{lem-pol-ana}
The convergence parameter $\alpha^{\ast}$ satisfies $0< \alpha^{\ast}\leq \alpha_{1}$. If $\alpha^{\ast} < \alpha_{1}$, then $\alpha^{\ast}$ is necessarily a zero point of $\hat{H}_{1}(\alpha, Z_{0}(\alpha))$.

\end{lemma}

\proof
From Lemma \ref{lem-asy-ana-1}-(i), we know that $\phi_{c-1}(\alpha)$ is analytic on $D_{\alpha}$ and thus the convergence parameter $\alpha^{\ast}\leq \alpha_{1}$.  For the case $\alpha^{\ast}< \alpha_{1}$,
we can deduce from Lemma~\ref{lem-asy-ana-1}-(i) that
$\alpha^{\ast}$ is either a
zero point of $\hat{H}_{1}(\alpha, Z_{0}(\alpha))$ or a point such that $|Z_{0}(\alpha^*)|\geq \frac{c\mu}{\lambda}$. In the following we prove
$|Z_{0}(\alpha)|< \frac{c\mu}{\lambda}$ for $\alpha\in (0, \alpha_{1})$.

For $\alpha \leq \alpha_{1}$,
we have
\[
Z_{0}(\alpha)= Z_{+}(\alpha)=\frac{-\alpha r+\lambda+c\mu - \sqrt{(-\alpha r+\lambda+c\mu)^{2}-4c\lambda\mu}}{2\lambda},
\]
which is a strictly increasing function of $\alpha$.
Thus, for any $\alpha\in (0, \alpha_{1})$, we have
\begin{equation}\label{equ-z}
 1=Z_{0}(0) < Z_{0}(\alpha)<  Z_{0}(\alpha_{1})= \sqrt{\frac{c\mu}{\lambda}}<\frac{c\mu}{\lambda}.
\end{equation}
\BOX

In order to perform the subsequent asymptotic arguments using the technique of complex analysis, we need to make some assumptions, which are collected as follows.

\begin{assumption}\label{ass-1}
(i) The function $\hat{H}_{1}(\alpha, Z_{0}(\alpha))$ has at most one real zero point in $(0, \alpha_{1}]$, denoted  by $\tilde{\alpha}$ if such a zero exists.

(ii) The zero point $\tilde{\alpha}$ satisfies  $H_{2}( Z_{0}(\tilde{\alpha}))\psi(Z_{0}(\tilde{\alpha}))+ \hat{H}_{0}(\tilde{\alpha}, Z_{0}(\tilde{\alpha}))\neq 0$.

(iii) The unique $\tilde{\alpha}$ is a multiple zero of $k$ times for function $\hat{H}_{1}(\alpha, Z_{0}(\alpha))$, where $k \geq 1$ is an integer.
\end{assumption}

\begin{remark}\label{rem-ass}
(i) For any set of model parameters $c$, $\lambda$ and $\mu$, (i), (ii) and (iii) of Assumption \ref{ass-1} can be easily checked numerically.

(ii) In many cases, these assumptions are not necessary. For example, if $c\mu> \lambda(r+1)$, we can derive from the expression of $\hat{H}_{1}(\alpha, Z_{0}(\alpha))$ that the unique zero point $\tilde{\alpha}$ must be a simple zero point. In this case, (iii) of Assumption \ref{ass-1} is redundant. Moreover, we will show in Section~\ref{sec:6} that all (i), (ii) and (iii) of Assumption \ref{ass-1} are redundant for the special cases $c=1$ and $c=2$. Actually, our extensive numerical calculations (for many sets of $\lambda$, $\mu$ and $r$ values) suggest that all (i), (ii) and (iii) are redundant for a general case, but a rigorous proof is still not available at this moment.
\end{remark}

The next lemma, which follows from Lemma~\ref{lem-asy-ana-1}-(i), provides more details about the convergence parameter $\alpha^{\ast}$.
\begin{lemma}\label{lem-sin-ana}
 Suppose that (i) and (ii) of  Assumption \ref{ass-1} hold. Then

(i) if the zero point $\tilde{\alpha}$ exists and $\tilde{\alpha}< \alpha_{1}$, we have  $\alpha^{\ast}=\tilde{\alpha}$,

(ii) if the zero point $\tilde{\alpha}$ exists and $\tilde{\alpha}= \alpha_{1}$, we have $\alpha^{\ast}=\tilde{\alpha}=\alpha_{1}$,


(iii) if $\hat{H}_{1}(\alpha, Z_{0}(\alpha))$ has no real zero points in $(0, \alpha_{1}]$, we have $\alpha^{\ast}=\alpha_{1}$.
\end{lemma}

Based on the above analysis, we can provide the following tail asymptotic properties for $\phi_{c-1}(\alpha)$ and $\psi(z)$, which are the key for
characterizing exact tail asymptotics in the stationary distribution of the fluid queue.
\begin{theorem}\label{the-tail-asy-Pi1}
 Suppose that (i) and (ii) of  Assumption \ref{ass-1} hold. For the function $\phi_{c-1}(\alpha)$, a total of three types of asymptotics exists as
$\alpha$ approaches to $\alpha^{\ast}$, based on the detailed property of $\alpha^{\ast}$ stated in Lemma~\ref{lem-sin-ana}.

Case (i) If (i) of Lemma \ref{lem-sin-ana} and (iii) of  Assumption \ref{ass-1} hold, then
\[
   \lim_{\alpha\rightarrow \alpha^{\ast}}(\alpha^{\ast}-\alpha)^{k}\phi_{c-1}(\alpha)= c_{1},
\]
where
\[
c_{1}=\frac{H_{2}(Z_{0}(\alpha^{\ast}))\psi(Z_{0}(\alpha^{\ast}))+ \hat{H}_{0}(\alpha^{\ast}, Z_{0}(\alpha^{\ast}))}
 {\hat{H}^{(k)}_{1}(\alpha^{\ast}, Z_{0}(\alpha^{\ast}))},
\]
and $\hat{H}^{(k)}_{1}(\alpha^{\ast}, Z_{0}(\alpha^{\ast}))$ is the $k$th derivative of $\alpha^{\ast}$.


Case (ii) If (ii) of Lemma \ref{lem-sin-ana} holds, then
\[
   \lim_{\alpha\rightarrow \alpha^{\ast}}\sqrt{\alpha^{\ast}-\alpha}\cdot\phi_{c-1}(\alpha)=c_{2},
\]
where
\[
c_{2}=\frac{2\lambda [H_{2}( Z_{0}(\alpha^{\ast}))\psi(Z_{0}(\alpha^{\ast}))+ \hat{H}_{0}(\alpha^{\ast}, Z_{0}(\alpha^{\ast}))]}{\frac{\partial \hat{H}_{1}(\alpha^{\ast}, Z_{0}(\alpha^{\ast}))}{\partial Z_{0}(\alpha^{\ast})}\cdot
\sqrt{\alpha^{\ast}-\alpha}}.
\]

Case (iii) If (iii) of Lemma \ref{lem-sin-ana} holds, then
\[
   \lim_{\alpha\rightarrow \alpha^{\ast}}\sqrt{\alpha^{\ast}-\alpha}\cdot\phi_{c-1}'(\alpha)=c_{3},
\]
where
\[
c_{3}=\frac{\partial L(\alpha, z)}{\partial z}|_{(\alpha^{\ast}, Z_{0}(\alpha^{\ast}))} \frac{\sqrt{\alpha_{2}-\alpha_{1}}}{2\lambda},
\]
and $L(\alpha, z)=-\frac{H_{2}(\alpha, z)\psi(z)+ \hat{H}_{0}(\alpha, z)}{ \hat{H}_{1}(\alpha, z)}.$
\end{theorem}

\proof
(i) In this case, $\alpha^{\ast}= \tilde{\alpha}$ is a multiple zero, with degree $k$, of  $\hat{H}_{1}(\alpha, Z_{0}(\alpha))$.
From (\ref{equ-asy-ana-1}), 
we have
\[
 (\tilde{\alpha}-\alpha)^{k}\phi_{c-1}(\alpha)=-\frac{H_{2}(Z_{0}(\alpha))\psi(Z_{0}(\alpha))+ \hat{H}_{0}(\alpha, Z_{0}(\alpha))}
 {\hat{H}_{1}(\alpha, Z_{0}(\alpha))/(\tilde{\alpha}-\alpha)^{k}}.
\]
It follows that
\begin{equation}\label{equ-alpha}
   \lim_{\alpha\rightarrow \tilde{\alpha}}(\tilde{\alpha}-\alpha)^{k}\phi_{c-1}(\alpha)=
 \frac{H_{2}( Z_{0}(\tilde{\alpha}))\psi(Z_{0}(\tilde{\alpha}))+ \hat{H}_{0}(\tilde{\alpha}, Z_{0}(\tilde{\alpha}))}
 {\hat{H}^{(k)}_{1}(\tilde{\alpha}, Z_{0}(\tilde{\alpha}))}=c_{1}.
\end{equation}
Moreover, as stated in Remark \ref{rem-ass}, we can obtain that  $c_{1}\neq 0$. Similarly, we also haven $c_{2}, c_{3} \neq 0$ in the following proof.

(ii) In this case, $\alpha^{\ast}= \tilde{\alpha}= \alpha_{1}$, which implies that
$\alpha_{1}$ is not only a zero point of $\Delta(\alpha)$ but also the zero point of $\hat{H}_{1}(\alpha, Z_{0}(\alpha))$. 
 Suppose that $\alpha_{1}$ is a zero of degree $m \geq 2$.  Then, we have
\[
 \lambda A_{c-2}'(\alpha_{1})=\frac{3}{2(c-1)\mu}r+\frac{1}{(c-1)\mu}> 0,
\]
which  conflicts the fact that $\lambda A_{c-2}'(\alpha_{1})< 0$. Hence, $\alpha_{1}$ is a simple zero point of $\hat{H}_{1}(\alpha, Z_{0}(\alpha))$.
%


Thus, we have
\begin{eqnarray*}
 \lim_{\alpha\rightarrow \alpha^{\ast}} \sqrt{\alpha^{\ast}-\alpha} \cdot\phi_{c-1}(\alpha) &=&  \lim_{\alpha\rightarrow \alpha^{\ast}}  -\frac{H_{2}( Z_{0}(\alpha))\psi(Z_{0}(\alpha))+ \hat{H}_{0}(\alpha, Z_{0}(\alpha))}{\hat{H}_{1}(\alpha, Z_{0}(\alpha)/\sqrt{\alpha^{\ast}-\alpha} } \\
    &=&  \lim_{\alpha\rightarrow \alpha^{\ast}} \frac{H_{2}(Z_{0}(\alpha))\psi(Z_{0}(\alpha))+ \hat{H}_{0}(\alpha, Z_{0}(\alpha))}{\sqrt{\alpha^{\ast}-\alpha}\cdot[\frac{\partial \hat{H}_{1}(\alpha^{\ast}, Z_{0}(\alpha^{\ast})}{\partial \alpha^{\ast}} + Z_{0}'(\alpha^{\ast})\cdot \frac{\partial \hat{H}_{1}(\alpha^{\ast}, Z_{0}(\alpha^{\ast}))}{\partial Z_{0}(\alpha^{\ast})}]},
\end{eqnarray*}
where  
\[
\lim_{\alpha\rightarrow \alpha^{\ast}}\sqrt{\alpha^{\ast}-\alpha}\cdot\frac{\partial \hat{H}_{1}(\alpha^{\ast}, Z_{0}(\alpha^{\ast})}{\partial \alpha^{\ast}}=0
\]
and
\begin{eqnarray*}
 \lim_{\alpha\rightarrow \alpha^{\ast}}\sqrt{\alpha^{\ast}-\alpha}Z_{0}'(\alpha^{\ast})
 &=& \lim_{\alpha\rightarrow \alpha^{\ast}}\frac{Z_{0}(\alpha^{\ast})-Z_{0}(\alpha)}{\sqrt{\alpha^{\ast}-\alpha}}\\
  &=& \lim_{\alpha\rightarrow \alpha^{\ast}}\left[\frac{b(\alpha^{\ast})-b(\alpha)}{2\lambda\sqrt{\alpha^{\ast}-\alpha}}+ \frac{\sqrt{(\alpha-\alpha^{\ast})(\alpha-\alpha_{2})}}{2\lambda\sqrt{\alpha^{\ast}-\alpha}}\right]\\
  &=& \frac{\sqrt{\alpha_{2}-\alpha^{\ast}}}{2\lambda}.
\end{eqnarray*}

It follows that
\[
 \lim_{\alpha\rightarrow \alpha^{\ast}}\sqrt{\alpha^{\ast}-\alpha}\cdot\phi_{c-1}(\alpha)=\frac{2\lambda [H_{2}( Z_{0}(\alpha^{\ast}))\psi(Z_{0}(\alpha^{\ast}))+ \hat{H}_{0}(\alpha^{\ast}, Z_{0}(\alpha^{\ast}))]}{\frac{\partial \hat{H}_{1}(\alpha^{\ast}, Z_{0}(\alpha^{\ast}))}{\partial Z_{0}(\alpha^{\ast})}\cdot
\sqrt{\alpha_{2}-\alpha^{\ast}}}=c_{2}.
\]

(iii) In this case, $\alpha^{\ast}= \alpha_{1}$. 
Let
\[
L(\alpha, z)=-\frac{H_{2}(z)\psi(z)+ \hat{H}_{0}(\alpha, z)}{ \hat{H}_{1}(\alpha, z)}.
\]
From (\ref{equ-asy-ana-1}), we have
\begin{eqnarray*}
  \phi_{c-1}'(\alpha) &=& \frac{\partial L(\alpha, z)}{\partial \alpha}+ \frac{\partial L(\alpha, z)}{\partial z} \cdot  Z_{0}'(\alpha).
\end{eqnarray*}
It follows that
\begin{eqnarray*}
 \lim_{\alpha\rightarrow \alpha^{\ast}}\sqrt{\alpha^{\ast}-\alpha}\cdot\phi_{c-1}'(\alpha) &=& \lim_{\alpha\rightarrow \alpha^{\ast}}\sqrt{\alpha^{\ast}-\alpha}\cdot\left[\frac{\partial L(\alpha, z)}{\partial \alpha}+ \frac{\partial L(\alpha, z)}{\partial z}\cdot Z_{0}'(\alpha)\right]\\
   &=&  \lim_{\alpha\rightarrow \alpha^{\ast}}\frac{\partial L(\alpha, z)}{\partial z}\cdot \sqrt{\alpha^{\ast}-\alpha}\cdot Z_{0}'(\alpha)\\
   &=&  \frac{\partial L(\alpha, z)}{\partial z}|_{(\alpha^{\ast}, Z_{0}(\alpha^{\ast}))} \frac{\sqrt{\alpha_{2}-\alpha_{1}}}{2\lambda}= c_{3}.
\end{eqnarray*}
\BOX

The asymptotic property for $\psi(z)$ can be stated as follows.
\begin{theorem}\label{the-tail-asy-Pi2}
For the function $\psi(z)$, we have the following asymptotic property as
$z$ approaches to $ \tilde{z}=\frac{c\mu}{\lambda}$:
\[
    \lim_{z\rightarrow \tilde{z}}(\tilde{z}-z)\psi(z)= d_{\tilde{z}},
\]
where
\[
 d_{\tilde{z}}= \frac{\hat{H}_{1}(\alpha(\tilde{z}), \tilde{z}) \phi_{c-1}(\alpha(\tilde{z}))+\hat{H}_{0}(\alpha(\tilde{z}), \tilde{z})}{\lambda(\tilde{z}-1)}.
\]
\end{theorem}
\proof
From (\ref{equ-asy-ana-2}), we have
\begin{eqnarray*}
  \lim_{z\rightarrow \tilde{z}}(\tilde{z}-z)\psi(z) &=&  \lim_{z\rightarrow \tilde{z}} \frac{\hat{H}_{1}(\alpha(z), z) \phi_{c-1}(\alpha(z))+\hat{H}_{0}(\alpha(z), z)}{\lambda(z-1)}\\
   &=& \frac{\hat{H}_{1}(\alpha(\tilde{z}), \tilde{z}) \phi_{c-1}(\alpha(\tilde{z}))+\hat{H}_{0}(\alpha(\tilde{z}), \tilde{z})}{\lambda(\tilde{z}-1)}.
\end{eqnarray*}
Moreover, we can calculate that  $\hat{H}_{1}(\alpha(\tilde{z}), \tilde{z})> 0$, which implies that $d_{\tilde{z}}> 0$.
\BOX

\section{Exact tail asymptotics for $\Pi_{i}(x)$ and $\Pi(x)$} \label{sec:5}

Lemma \ref{lem-tail-asy-2} and Lemma \ref{lem-tail-asy-Pi} specify exact tail asymptotic properties for the density function $\pi_{c-1}(x)$ and boundary probabilities $\Pi_{i}(0)$, respectively, which are direct consequences of the detailed asymptotic behavior of $\phi_{c-1}(\alpha)$ and $\psi(z)$, and the Tauberian-like theorem, given in Appendix~\ref{app:A}. Moreover, the tail asympotics for the joint probability $\Pi_{i}(x)$, the density function $\pi_{i}(x)$, the marginal distribution $\Pi(x)=\sum_{i=0}^{\infty}\Pi_{i}(x)$,  and the density function $\pi(x)= \frac{d\Pi(x)}{dx}$ for $x> 0$ are also provided in this section.

\begin{lemma}\label{lem-tail-asy-2}
 Suppose that (i) and (ii) of  Assumption \ref{ass-1} hold. For the density function $\pi_{c-1}(x)$ of the fluid queue, we have the following tail asymptotic properties for large enough
$x$.

Case (i) If (i) of Lemma \ref{lem-sin-ana} and (iii) of  Assumption \ref{ass-1} hold, then
\[
   \pi_{c-1}(x)\sim C_{1}e^{-\alpha^{\ast} x}x^{k-1}.
\]

Case (ii) If (ii) of Lemma \ref{lem-sin-ana} holds, then
\[
   \pi_{c-1}(x) \sim C_{2} e^{-\alpha^{\ast} x}x^{-\frac{1}{2}}.
\]

Case (iii) If (iii) of Lemma \ref{lem-sin-ana} holds, then
\[
  \pi_{c-1}(x) \sim C_{3} e^{-\alpha^{\ast} x}x^{-\frac{3}{2}},
\]
where 
$C_{1}=\frac{c_{1}}{\Gamma(k)}$, $C_{2}=\frac{c_{2}}{\sqrt{\pi}}$, $C_{3}=\frac{-c_{3}}{2\sqrt{\pi}}$  and $c_{i}, i=1, 2, 3$ are defined in Theorem \ref{the-tail-asy-Pi1}.
\end{lemma}
\begin{lemma}\label{lem-tail-asy-Pi}
For the boundary probabilities $\Pi_{i}(0)$ of the fluid queue, we have the following tail asymptotic properties for large enough
$i$.
 \[
      \Pi_{i}(0)\sim d_{\tilde{z}}\cdot \Big (\frac{1}{\tilde{z}} \Big )^{i+1},
 \]
where $\tilde{z}=\frac{c\mu}{\lambda}$ and $d_{\tilde{z}}$ is defined in Theorem \ref{the-tail-asy-Pi2}.
\end{lemma}
Now we provide details for the exact tail asymptotic characterization in the (general) joint probabilities $\Pi_{i}(x)$ for any
 $i\geq c-1$.

\begin{theorem}\label{the-tail-asy-3}
 Suppose that (i) and (ii) of  Assumption \ref{ass-1} hold. For the joint probabilities $\Pi_{i}(x)$ of the fluid queue, then we have the following tail asymptotic properties for any $i \geq c-1$ and large enough $x$.

Case (i) If (i) of Lemma \ref{lem-sin-ana} and (iii) of  Assumption \ref{ass-1} hold, then
\begin{equation}\label{equ-joint-probab}
   \Pi_{i}(x)- \xi_{c} \Big (\frac{\lambda}{c\mu} \Big )^{i-c} \sim -\frac{C_{1}}{\alpha^{\ast}}e^{-\alpha^{\ast}x}x^{k-1} \Big (\frac{1}{z^{\ast}} \Big )^{i-c},
\end{equation}
and
\[
  \pi_{i}(x) \sim C_{1}e^{-\alpha^{\ast}x}x^{k-1} \Big (\frac{1}{z^{\ast}} \Big )^{i-c}.
\]

Case (ii) If (ii) of Lemma \ref{lem-sin-ana}  holds, then
\[
   \Pi_{i}(x)-  \xi_{c} \Big (\frac{\lambda}{c\mu} \Big )^{i-c} \sim -\frac{C_{2}}{\alpha^{\ast}} e^{-\alpha^{\ast}x}x^{-\frac{1}{2}} \Big (\frac{1}{z^{\ast}} \Big )^{i-c},
\]
and
\[
  \pi_{i}(x) \sim C_{2}e^{-\alpha^{\ast}x}x^{-\frac{1}{2}} \Big (\frac{1}{z^{\ast}} \Big )^{i-c}.
\]

Case (iii) If (iii) of Lemma \ref{lem-sin-ana}  holds, then
\[
   \Pi_{i}(x)-  \xi_{c} \Big (\frac{\lambda}{c\mu} \Big )^{i-c} \sim -\frac{C_{3}}{\alpha^{\ast}} e^{-\alpha^{\ast}x}x^{-\frac{3}{2}} \Big (\frac{1}{z^{\ast}} \Big )^{i-c},
\]
and
\[
  \pi_{i}(x) \sim C_{3}e^{-\alpha^{\ast}x}x^{-\frac{3}{2}} \Big (\frac{1}{z^{\ast}} \Big )^{i-c},
\]
where 
 $z^{\ast}=Z_{0}(\alpha^{\ast})$ and  $C_{1}$, $C_{2}$, $C_{3}$ are  defined in Lemma \ref{lem-tail-asy-2}.
\end{theorem}

\proof
We only prove (i), since (ii) and (iii) can be similarly proved. 
For  $i=c-1$, we have
\[
 \lim_{x\rightarrow\infty} \frac{C_{1}e^{-\alpha^{\ast} x}x^{k-1}}{\xi_{c-1}-\Pi_{c-1}(x)}
    = \lim_{x\rightarrow\infty} \frac{\alpha^{\ast}C_{1}e^{-\alpha^{\ast} x}x^{k-1}-(k-1)C_{1}e^{-\alpha^{\ast} x}x^{k-2}}{\pi_{c-1}(x)}
    =  \alpha^{\ast},
\]
where the first equality follows from the L'Hospital's rule and the second equality follows from Lemma \ref{lem-tail-asy-2}. Hence, we have $\Pi_{c-1}(x)-\xi_{c-1} \sim -\frac{C_{1}}{\alpha^{\ast}}e^{-\alpha^{\ast}x}x^{k-1}$ as $x\rightarrow \infty$.

Now, suppose that (\ref{equ-joint-probab}) is true for any $i=m > c-1$. Thus,
for
$i=m+1$ it follows from (\ref{equ-joi-dis-3}) that
\[
  c\mu \Pi_{m+1}(x)=
  -\lambda\Pi_{m-1}(x)+(\lambda+c\mu)\Pi_{m}(x)+r \pi_{m}(x),
\]
which leads to
\begin{eqnarray*}
  &&\lim_{x\rightarrow\infty}\frac{ \Pi_{m+1}(x)-\xi_{c}(\frac{\lambda}{c\mu})^{m+1-c}}{\frac{C_{1}}{\alpha^{\ast}}e^{-\alpha^{\ast}x}x^{k-1}} \\
  &=&\lim_{x\rightarrow\infty} \left [-\frac{\lambda}{c\mu}\cdot\frac{\Pi_{m-1}(x)-\xi_{c}(\frac{\lambda}{c\mu})^{m-c-1}}{\frac{C_{1}}{\alpha^{\ast}}
  e^{-\alpha^{\ast}x}x^{k-1}}
  +\frac{\lambda+c\mu}{c\mu}\cdot\frac{ \Pi_{m}(x)-\xi_{c}(\frac{\lambda}{c\mu})^{m-c}}{\frac{C_{1}}{\alpha^{\ast}}e^{-\alpha^{\ast}x}x^{k-1}}+
  \frac{r\alpha^{\ast}}{c\mu}\cdot\frac{ \pi_{m}(x)}{C_{1}e^{-\alpha^{\ast}x}x^{k-1}} \right ]\\
  &=&  -(\frac{1}{z^{\ast}})^{m-c} \left [-\frac{\lambda}{c\mu}z^{\ast}+\frac{\lambda+ c\mu}{c\mu}- \frac{r \alpha^{\ast}}{c\mu } \right ] \\
  &=&  -(\frac{1}{z^{\ast}})^{m+1-c},
\end{eqnarray*}
where the last equation follows from the fact that $H(\alpha^{\ast}, z^{\ast})=0$ and $z^{\ast}=Z_{0}(\alpha^{\ast})$. This completes the proof.

\BOX
\begin{remark}\label{rem-c}
According to (\ref{equ-phi-rel}), we can derive tail asymptotic properties of $\phi_{c-2}(\alpha)$ from $\phi_{c-1}(\alpha)$, and
thus tail asymptotic properties for $\Pi_{c-2}(x)$. Similarly,  a relationship can be established between
$\phi_{i}(\alpha)$ and $\phi_{i-1}(\alpha)$ for any $1 \leq i\leq c-2 $, and thus tail asymptotic properties for the joint probability $\Pi_{i}(x)$
can be obtained for any $0\leq i \leq c-1$.
\end{remark}

In the following theorem, we provide exact tail asymptotics for the marginal probabilities
 $\Pi(x)$.
\begin{theorem}\label{the-tail-asy-4}
 Suppose that (i) and (ii) of  Assumption \ref{ass-1} hold. For the marginal probabilities $\Pi(x)$ of the fluid queue, we have the following tail asymptotic properties: 

Case (i) If (i) of Lemma \ref{lem-sin-ana} and (iii) of  Assumption \ref{ass-1} hold, then
\begin{equation*}
   \Pi(x)- 1 \sim -\frac{\tilde{C}_{1}}{\alpha^{\ast}}e^{-\alpha^{\ast}x}x^{k-1},
\end{equation*}
and
\[
  \pi(x) \sim \tilde{ C}_{1}e^{-\alpha^{\ast}x}x^{k-1}.
\]

Case (ii) If (ii) of Lemma \ref{lem-sin-ana}  holds, then
\[
    \Pi(x)- 1 \sim -\frac{\tilde{C}_{2}}{\alpha^{\ast}}e^{-\alpha^{\ast}x}x^{-\frac{1}{2}},
\]
and
\[
  \pi(x) \sim e^{-\alpha^{\ast}x}x^{-\frac{1}{2}}.
\]

Case (iii) If (iii) of Lemma \ref{lem-sin-ana}  holds, then
\[
   \Pi(x)-  1 \sim -\frac{\tilde{C}_{3}}{\alpha^{\ast}}e^{-\alpha^{\ast}x}x^{-\frac{3}{2}},
\]
and
\[
  \pi(x) \sim \tilde{C}_{3}e^{-\alpha^{\ast}x}x^{-\frac{3}{2}},
\]
where 
  $\tilde{C}_{i}=\left[\frac{\hat{H}_{1}(\alpha^{\ast}, 1)}{ H(\alpha^{\ast}, 1)}+\sum_{k=0}^{c-2}A_{k}(\alpha^{\ast}) A_{k+1}(\alpha^{\ast})\cdots A_{c-2}(\alpha^{\ast})\right]C_{i}$, $i=1, 2, 3$, $C_{1}$, $C_{2}$ and $C_{3}$ are  defined in Lemma~\ref{lem-tail-asy-2}, and $A_{i}(\alpha)$ is defined in Theorem~\ref{the-phi-rel}.
\end{theorem}

\proof
Let $z=1$. It follows from (\ref{equ-funde}) that
\[
 H(\alpha, 1)\psi(\alpha, 1)= \hat{H}_{1}(\alpha, 1)\phi_{c-1}(\alpha)+  H_{2}( 1)\psi(1)+ \hat{H}_{0}(\alpha, 1).
\]
Thus, we get
\begin{equation}\label{equ-mar-dis}
   H(\alpha, 1)\int_{0}^{\infty}\sum_{i=c-1}^{\infty} \pi_{i}(x) e^{\alpha x}dx=\hat{H}_{1}(\alpha, 1)\int_{0}^{\infty}\pi_{c-1}(x) e^{\alpha x}dx+ \hat{H}_{0}(\alpha, 1),
\end{equation}
since $H_2(1)=0$.
From (\ref{equ-mar-dis}), we have
\[
  \int_{0}^{\infty} \sum_{i=c-1}^{\infty} \pi_{i}(x) e^{\alpha x}dx=\frac{\hat{H}_{1}(\alpha, 1)}{ H(\alpha, 1)}\int_{0}^{\infty}\pi_{c-1}(x) e^{\alpha x}dx+\frac{\hat{H}_{0}(\alpha, 1)}{H(\alpha, 1)}.
\]
Now, we prove
\begin{equation}\label{equ-exc-orde}
  \int_{0}^{\infty} \sum_{i=0}^{\infty} \pi_{i}(x) e^{\alpha x}dx
  =\int_{0}^{\infty} \pi(x) e^{\alpha x}dx,
\end{equation}
where $\pi_{i}(x)=\frac{\partial\Pi_{i}(x)}{\partial x}$ and $\pi(x)=\frac{d\Pi(x)}{dx}$ for any $x> 0$.

For any fixed $x$, we can obtain
\[
  \sum_{i=0}^{\infty} \Pi_{i}(x)= P\{X< x\}\leq 1,
\]
which implies that $\sum_{i=0}^{\infty} \Pi_{i}(x)$ is convergent for any $x$. From (\ref{equ-joi-dis-3}), we have for $i\geq c$,
\[
  \pi_{i}(x)=\frac{\lambda}{r}\Pi_{i-1}(x)- \frac{\lambda+ c\mu}{r}\Pi_{i}(x)+ \frac{c\mu}{r} \Pi_{i+1}(x)\leq \frac{\lambda}{r} \xi_{i-1}+\frac{c\mu}{r} \xi_{i+1}.
\]
Since
\[
\sum_{i=c}^{\infty} \frac{\lambda}{r} \xi_{i-1}+ \sum_{i=c}^{\infty} \frac{c\mu}{r} \xi_{i+1} < \frac{\lambda+ c\mu}{r}< \infty,
\]
according to Weierstrass criterion,  we can obtain that $\sum_{i=0}^{\infty} \pi_{i}(x)$ is convergent uniformly in $x$.
Thus, we can get equation (\ref{equ-exc-orde}). From (\ref{equ-mar-dis}), we have
\[
 \int_{0}^{\infty} \pi(x) e^{\alpha x}dx=\frac{\hat{H}_{1}(\alpha, 1)}{ H(\alpha, 1)}\phi_{c-1}(\alpha)+\frac{\hat{H}_{0}(\alpha, 1)}{H(\alpha, 1)}+\sum_{i=0}^{c-2}\phi_{i}(\alpha).
\]
From (\ref{equ-phi-rel}), we can establish the relationship between $\phi_{i}(\alpha)$ and $\phi_{c-1}(\alpha)$ for any $0< i\leq c-2$, and thus
\[
 \sum_{i=0}^{c-2}\phi_{i}(\alpha)=\left(\sum_{k=0}^{c-2}A_{k}(\alpha)A_{k+1}(\alpha)\cdots A_{c-2}(\alpha)\right)\phi_{c-1}(\alpha)+ H_{c-1}(\alpha).
\]
Here $H_{c-1}(\alpha)$ is an analytic function about $\alpha$, which can be determined explicitly by (\ref{equ-phi-rel}). Hence, according to the Tauberian-like theorem and the asymptotic behavior of $\phi_{c-1}(\alpha)$, we can obtain
the tail asymptotic properties of $\pi(x)$ and thus attain the tail asymptotic properties of $\Pi(x)$.
\BOX

\section{Special cases} \label{sec:6}

In this section, we consider two important special cases: $c=1$ and $c=2$, for which exact asymptotic properties for the stationary distribution can be obtained without Assumption~\ref{ass-1}. The analysis of these two cases is feasible. However, the arguments for the cases $c\geq 3$ are rather complex since
the expression of $\hat{H}_{1}(\alpha, Z_{0}(\alpha))$ is intractable for any $c\geq 3$.

\subsection{Fluid queue driven by $M/M/1$ queue}

In this case, the unique zero point of $\hat{H}_{1}(\alpha, Z_{0}(\alpha))$ can be obtained explicitly as follows.
\begin{lemma}\label{lem-asy-ana-6}
Let
\[
  \tilde{\alpha}=\frac{\mu}{r+1}-\lambda,
\]
then $\tilde{\alpha}$ is the only possible zero point of  $H_{1}(\alpha, Z_{0}(\alpha))$. Moreover, $\tilde{\alpha}$ must be a simple
zero point of $H_{1}(\alpha, Z_{0}(\alpha))$.
 \end{lemma}

\proof
We rationalize  $\hat{H}_{1}(\alpha, Z_{0}(\alpha))$  by
\begin{equation}\label{equ-zero}
  g(\alpha)=2a\hat{H}_{1}(\alpha, Z_{0}(\alpha))\hat{H}_{1}(\alpha, Z_{1}(\alpha)).
\end{equation}
 Then, it follows from  the definition of $\hat{H}_{1}(\alpha, z)$ in Theorem \ref{the-phi-rel} and (\ref{equ-Z})
that
\begin{eqnarray*}
  g(\alpha) &=& -2\lambda Z_{0}(\alpha)Z_{1}(\alpha)[(\mu -\alpha r - \alpha) Z_{0}(\alpha)- \mu][(\mu -\alpha r - \alpha) Z_{1}(\alpha)- \mu] \\
    &=& \frac{-2\alpha\mu^{2}}{\lambda}[( r +1)\alpha-\mu+ \lambda(r+1)].
\end{eqnarray*}
It is obvious that $\tilde{\alpha}=\frac{\mu}{r+1}-\lambda$ is the only possible zero point  of $H_{1}(\alpha, Z_{0}(\alpha))$ with a modulus greater than $0$.
\BOX

\begin{lemma}\label{lem-asy-ana-6}
The unique zero point $\tilde{\alpha}$ satisfies the following inequality
\[
  H_{2}( Z_{0}(\tilde{\alpha}))\psi(Z_{0}(\tilde{\alpha}))+ \hat{H}_{0}(\tilde{\alpha}, Z_{0}(\tilde{\alpha}))\neq 0.
\]
 \end{lemma}

\proof
From the initial condition $\Pi_{i}(0)= 0$  for any $i\geq 1$, we have
\[
    \psi(Z_{0}(\tilde{\alpha}))=\sum_{i=0}^{\infty}\Pi_{i}(0)Z_{0}^{i}(\tilde{\alpha})=\Pi_{0}(0).
\]
For $\tilde{\alpha}=\frac{\mu}{r+1}-\lambda$, we have $Z_{0}(\tilde{\alpha})=\min\{1+r, \frac{\mu}{\lambda(1+r)}\}$.
Thus, from the definitions of $\hat{H}_{0}(\alpha, z)$ and $H_{2}(z)$ in Theorem \ref{the-phi-rel}, and the fact that $\frac{\mu}{\mu-\lambda Z_{0}(\tilde{\alpha})}\neq 1$, 
we can obtain that
\[
 \frac{-\hat{H}_{0}(\tilde{\alpha}, Z_{0}(\tilde{\alpha}))}{ H_{2}(Z_{0}(\tilde{\alpha}))}=\frac{\mu\Pi_{0}(0)}{\mu-\lambda Z_{0}(\tilde{\alpha})}\neq \Pi_{0}(0), 
\]
which implies that $H_{2}(Z_{0}(\tilde{\alpha}))\psi(Z_{0}(\tilde{\alpha}))+ \hat{H}_{0}(\tilde{\alpha}, Z_{0}(\tilde{\alpha}))\neq 0$.
\BOX

Since the following  inequality always holds
\[
\tilde{\alpha}= \frac{\mu}{r+1}-\lambda\leq \alpha_{1}=\frac{(\sqrt{\mu}-\sqrt{\lambda})^{2}}{r},
\]
 we only have two tail asymptotic properties for the stationary distribution $\Pi_{i}(x)$ and the marginal distribution $\Pi(x)$ of the fluid queue in this case. Here we omit the details and only present the asymptotic property for the marginal distribution.

 \begin{theorem}\label{the-tail-asy-marginal-6}
For the marginal distribution $\Pi(x)$ of the fluid queue, we have the following tail asymptotic properties for large enough
$x$:

Case (i) If (i) of Lemma \ref{lem-sin-ana} holds (i.e. $\tilde{\alpha}= \frac{\mu}{r+1}-\lambda<\alpha_{1}$), then 
\[
  \Pi(x)- 1 \sim -\frac{(r+1)\tilde{\alpha}c_{1}}{r}e^{-\tilde{\alpha} x};
\]

Case (ii) If (ii) of Lemma \ref{lem-sin-ana} holds (i.e. $\tilde{\alpha}= \frac{\mu}{r+1}-\lambda=\alpha_{1}$), then
\[
  \Pi(x)- 1 \sim -\frac{(r+1)\tilde{\alpha}c_{2}}{r\sqrt{\pi}} e^{-\tilde{\alpha} x}x^{-\frac{1}{2}}.
\]
Here  $c_{1}$ and $c_{2}$ are defined in Theorem~\ref{the-tail-asy-Pi1}.
\end{theorem}

\subsection{Fluid queue driven by $M/M/2$ queue}

%
 In this case, we can obtain the following lemma.
 \begin{lemma}\label{lem-asy-ana-6-b}
The function $\hat{H}_{1}(\alpha, Z_{0}(\alpha))$ has at most one real zero point in $(0, \alpha_{1}]$. 
 Moreover, this unique zero point, denoted  by $\tilde{\alpha}$ if exists, must be
a simple zero point.
\end{lemma}

\proof
Let $g(\alpha)=0$, where $g(\alpha)$ is defined in (\ref{equ-zero}). We can obtain the following equation:
\[
 (r+1)\alpha^{3}+[3\lambda(r+1)+\mu r]\alpha^{2}+[3\lambda^{2}(r+1)+\mu\lambda r-\lambda\mu-\mu^{2}]\alpha+\lambda^{3}(r+1)-\lambda^{2}\mu-2\lambda\mu^{2}=0.
\]
Denote the left hand side of the above equation by $\tilde{g}(\alpha)$.  For any $\alpha> 0$, we have
\[
\tilde{g}''(\alpha)=6(r+1)\alpha+6\lambda(r+1)+2\mu r> 0,
\]
which implies that $\tilde{g}(\alpha)$ is a convex function for any $\alpha> 0$.

If $\tilde{g}(0)=\lambda^{3}(r+1)-\lambda^{2}\mu-2\lambda\mu^{2} < 0$,  we can derive that  $\tilde{g}(\alpha)=0$ has only one real solution on $(0, \infty)$, which implies that there exists at most one real solution on $(0, \alpha_{1}]$, denoted by  $\tilde{\alpha}$ if exists. Moreover, according to the property of convex functions, we have $\tilde{g}'(\tilde{\alpha})>0$, which implies
 that $\tilde{\alpha}$ is a simple zero point.

 If $\tilde{g}(0) \geq 0$, we have  $\tilde{g}'(0)=3\lambda^{2}(r+1)+\mu\lambda r-\lambda\mu-\mu^{2}> 0$ and thus
$\tilde{g}(\alpha)=0$ has no real solution on $(0, \infty)$.
\BOX
%


\begin{lemma}\label{lem-asy-ana-6}
The unique zero point  $\tilde{\alpha}$, if it exists, satisfies the following inequality
\[
  H_{2}( Z_{0}(\tilde{\alpha}))\psi(Z_{0}(\tilde{\alpha}))+ \hat{H}_{0}(\tilde{\alpha}, Z_{0}(\tilde{\alpha}))\neq 0.
\]
 \end{lemma}

\proof
From the definitions of $H_{2}(z)$ and $\hat{H}_{0}(\alpha, z)$ in Theorem \ref{the-phi-rel}, for $c=2$,  we have
\begin{eqnarray}\nonumber
  &&H_{2}(Z_{0}(\alpha))\psi(Z_{0}(\alpha))+ \hat{H}_{0}(\alpha, Z_{0}(\alpha)) \\
\nonumber &=& (\lambda Z_{0}(\alpha)-2\mu)(Z_{0}(\alpha)-1)\psi(Z_{0}(\alpha))+ \left[\mu Z_{0}(\alpha)^{2}-2\mu Z_{0}(\alpha)+\frac{\lambda \mu Z_{0}^{2}(\alpha)}{\alpha+\lambda}\right]\Pi_{1}(0)+\frac{\lambda \alpha Z_{0}^{2}(\alpha)}{\alpha+\lambda}\Pi_{0}(0)\\
 \label{equ-not-zero} &=&  Z_{0}(\alpha)^{2} \left[\lambda (Z_{0}(\alpha)-1)\Pi_{1}(0)+\frac{\alpha (\lambda\Pi_{0}(0)-\mu \Pi_{1}(0))}{\alpha+\lambda}\right],
\end{eqnarray}
where the second equality follows from the fact that $\psi(z)=\Pi_{1}(0)z$.

Actually, from (\ref{equ-joi-dis-1}), we can obtain that $\lambda\Pi_{0}(0)-\mu \Pi_{1}(0)\geq 0$. Moreover, from (\ref{equ-z}), we have $Z_{0}(\alpha)> 1$ for any $\alpha \in (0, \alpha_{1}]$.
Hence, from (\ref{equ-not-zero}) we can get
\[
  H_{2}( Z_{0}(\tilde{\alpha}))\psi(Z_{0}(\tilde{\alpha}))+ \hat{H}_{0}(\tilde{\alpha}, Z_{0}(\tilde{\alpha}))> 0.
\]
\BOX

From the above lemmas, we can obtain the following theorem.
\begin{theorem}\label{the-tail-asy-6}
For the marginal probabilities $\Pi(x)$ of the fluid queue, we have the following tail asymptotic properties: 

Case (i) If (i) of Lemma \ref{lem-sin-ana} holds, then
\begin{equation}\label{equ-mar-probab-2}
   \Pi(x)- 1 \sim -\frac{\tilde{C}_{1}}{\alpha^{\ast}}e^{-\alpha^{\ast}x};
\end{equation}

Case (ii) If (ii) of Lemma \ref{lem-sin-ana} holds, then
\[
    \Pi(x)- 1 \sim -\frac{\tilde{C}_{2}}{\alpha^{\ast}}e^{-\alpha^{\ast}x}x^{-\frac{1}{2}};
\]

Case (iii)   If (iii) of Lemma \ref{lem-sin-ana} holds, then 
\[
   \Pi(x)-  1 \sim -\frac{\tilde{C}_{3}}{\alpha^{\ast}}e^{-\alpha^{\ast}x}x^{-\frac{3}{2}}.
\]
Here $\tilde{C}_{i}$, $i=1, 2, 3$, are defined in~Theorem \ref{the-tail-asy-4}.
\end{theorem}

\begin{remark}\label{rem-c2}
  Compared with the case of  $c=1$, a new asymptotic behavior, Case (iii), appears in the case of $c=2$.  We now give an example to illustrate that this new asymptotics exists for the case of $c \geq 3$.  For example, let $c=3$, if we take $r=10,$ $\lambda=20$, $\mu=30$, we can obtain four zero points $4,$ $-67$ and $-15\pm 5i$ of $\hat{H}_{1}(\alpha, Z_{0}(\alpha))$. Thus, we have $\tilde{\alpha}=4 > \alpha_{1}=0.5$, which implies that Case (iii) holds.
\end{remark}

\section{Concluding remarks} \label{sec:7}

In this paper, we applied the kernel method to investigate exact tail asymptotic properties of the joint stationary probabilities and the marginal distribution of the fluid queue driven by an $M/M/c$ queue. Different from the model studied in \cite{LZ11} and \cite{LTZ13}, for which tail asymptotic properties are symmetric between the level and the phase, since both level and phase processes are discrete,  the tail asymptotics for $\phi_{c-1}(\alpha)$ and $\psi(z)$ are asymmetric in this paper, since the phase process is discrete and level process is continuous.

There exists a total of three different types of exact tail asymptotics for the stationary probabilities of the
fluid queue in this paper. However, we may see in \cite{GLR13} that the stationary probabilities of the fluid queue driven by a finite Markov chain is always exactly geometric, which corresponds to Case (i) of Theorem \ref{the-tail-asy-3}. This implies that the infinite phase causes new phenomena.

In Section 6, we showed that Case (iii) of Theorem \ref{the-tail-asy-3} does not appear in the case of $c=1$, but exists for the case of
$c\geq 2$. This implies that the asymptotic behaviour for $c \geq 2$ can be significantly different from that for the case of $c=1$.

From the arguments given in the paper, we have seen that Assumption \ref{ass-1} is redundant in the special cases $c=1$ and $c=2$. Based on our numerical calculations for a broad range selections of parameter values, we conjecture that Assumption \ref{ass-1} is redundant for all cases of $c\geq 3$.

\section*{Acknowledgements}
This research was supported in part by the National Natural Science
Foundation of China (grant 11571372, 11771452),  and a Discovery Grant by the Natural Sciences and
Engineering Research Council of Canada (NSERC).

\appendix
\newpage

\section{Tauberian-like theorem} \label{app:A}

Denote
\[
  \Delta_{1}(\phi, \varepsilon)=\{x: |x|\leq |x_{0}|+\varepsilon, |\arg(x-x_{0})|> \phi, \varepsilon>0, 0<\phi<\frac{\pi}{2}\}.
\]

Let $f_{n}$ be a sequence of numbers, with the generating function
\[
 f(x)=\sum_{n\geq 1} f_{n} x^{n}.
\]

\begin{lemma}\label{lem-tauber-like-1}
(Flajolet and Odlyzko 1990)  Assume that $f(x)$ is analytic in $\Delta_{1}(\phi, \varepsilon)$ expect at $x=x_{0}$ and
\[
  f(x)\sim K(x_{0}-x)^{s} \ \ \mbox{as}\ \   x\rightarrow x_{0} \ \ \mbox{in} \ \  \Delta_{1}(\phi, \varepsilon).
\]
Then as $n\rightarrow \infty$,
(i)  If $s \not\in \{0, 1, 2, \ldots\}$,
  \[
   f_{n}=\frac{K }{\Gamma(-s)}n^{-s-1}x_{0}^{-n},
  \]
 where $\Gamma(\cdot)$ is the Gamma function.

(ii) If $s$ is a non-negative integer, then
\[
 f_{n}= o(n^{-s-1}x_{0}^{-n}).
\]
\end{lemma}

For the continuous case,
let
\[
g(x)=\int_{0}^{\infty}e^{xt}f(t)dt.
\]
Denote
\[
  \Delta_{2}(\phi, \varepsilon)=\{x: \Re(x)\leq |x_{0}|+\varepsilon, x\neq x_{0}, \varepsilon>0, |\arg(x-x_{0})|> \phi\}.
\]
The following lemma has been shown in Theorem 2 in \cite{DDZ15}.

\begin{lemma}\label{the-tauber-like-2}
Assume that $g(x)$ satisfies the following conditions:

(i) The left-most singularity of $g(x)$ is $x_{0}$ with $x_{0} > 0$. Furthermore, we assume that
as $x\rightarrow x_{0}$,
\[
  g(x)\sim (x_{0}-x)^{-s}
\]
for some $s\in \mathbb{C}\backslash Z_{-}$.

(ii) $g(x)$ is analytic on $\Delta_{2}(\phi_{0}, \varepsilon)$  for some $\phi_{0}\in (0, \frac{\pi}{2}]$.

(iii) $g(x)$ is bounded on $\Delta_{2}(\phi_{1}, \varepsilon)$ for some $\phi_{1}> 0$.

Then, as $t\rightarrow \infty$,
\[
  f(t)\sim e^{-x_{0} t} \frac{t^{s-1}}{\Gamma(s)},
\]
where $\Gamma(\cdot)$ is the Gamma function.
\end{lemma}
\end{document}